\input amstex
\documentstyle{amsppt}

\NoRunningHeads
\mag=1200
\hsize=28.5 pc
\vsize=41 pc
\hcorrection{2mm}
\topmatter
\title
On finiteness of the number of boundary slopes of immersed
surfaces in 3-manifolds
\endtitle
\author
Joel Hass and Shicheng Wang and Qing Zhou
\endauthor
\address The University of California at Davis, CA, 95616, USA \endaddress
\address Peking University, Beijing 100871, China \endaddress
\address East China Normal University, Shanghai, 200062, China \endaddress 

\abstract
For any hyperbolic 3-manifold $M$ with totally geodesic boundary,
there are finitely many boundary slopes for essential immersed surfaces of
a given genus. There is a uniform bound for the number of such boundary
slopes if the genus of $\partial M$ or the volume of $M$ is bounded above.
When the volume is bounded above, then area of $\partial M$ is bounded
above and the length of closed geodesic on $\partial M$
is bounded below.
\endabstract

\thanks
This paper grew out of work begun
while the first two authors were visiting the Mathematical
Sciences Research Institute in Berkeley in 1996-97. Research at MSRI
is supported in part by NSF grant DMS-9022140. The first author was
partially supported by NSF grant DMS-9704286. The
second and third authors were partially supported by  Outstanging Younth
Fellowship of NSFC \endthanks

\endtopmatter
\document

\normalbaselineskip=1.44\normalbaselineskip
\normalbaselines

We say that a proper immersion of a surface $F$ into $M$ is an
{\it essential surface} if it is
incompressible and $\partial$-incompressible,
meaning that the immersion induces
an injection of the fundamental group and relative
fundamental group. Let $c$ be an essential simple loop on
the boundary $\partial M$ of a compact 3-manifold $M$.
If there is a proper immersion of an essential surface $F$
into $M$ such that each component of $\partial F$ is homotopic
to a multiple of $c$, we call $c$ a {\it boundary slope} of $F$. 

We are interested in the following two questions: 

\proclaim{Questions}

\noindent
(1) Given a compact 3-manifold $M$ and a genus $g$, are there
finitely many boundary slopes for immersed essential surfaces with
genus at most $g$? 

\noindent
(2) Under what conditions is there a bound for the number
of boundary slopes in (1) which is independent of the 3-manifold? 

\endproclaim

Many results in these directions have been obtained for various classes
of 3-manifolds: 

\noindent
(1) If $\partial M$ is a torus and the surfaces are embedded,
Hatcher [H] showed that there are only finitely many boundary slopes,
without any genus restriction.

\noindent
(2) When the surfaces are embedded punctured spheres or tori,
explicit bounds are known on the number of boundary slopes.
These bounds are based on highly developed combinatorial methods
in knot theory and the theory of representations of knot groups.
See the survey papers [Go], [Lu] and [Sh]. 

\noindent
(3) When $\partial M$ is a torus and the surfaces are immersed,
a positive answer to Question (1) has been obtained recently in [HRW].
When $M$ is hyperbolic, minimal surface theory is used to derive these bounds.
For fixed genus $g$, these turn out to be quadratic functions of $g$,
independent of $M$. See also recent work of Agol [Agol]. 

\noindent
(4) If $M$ is an irreducible, $\partial$-irreducible, acylindrical,
atoroidal 3-manifold and the surfaces are embedded,
Scharlemann and Wu [SW] gave a positive answer to Question 1
using combinatorial arguments. 

\noindent
(5) Suppose $\partial M$ is a torus and the surfaces are immersed.
Baker has given examples to show that the bounded genus assumption cannot
be dropped. Oertel, using branched surface theory,
has found manifolds in which every slope is realized by the boundary
of an immersed essential surface [Oe].

In this note we give a positive answer to Question (1),
which extends (3) to the case where $\partial M$ can contain
high genus components and generalizes (4) from embedded to immersed surfaces. 

\proclaim{Theorem 1}
Suppose $M$ is
$\partial$-irreducible, acylindrical, atoroidal 3-manifold.
Then for any $g$, there are only finitely many $\partial$-slopes
for essential surfaces of genus $g$. \endproclaim 

Next we consider the question of obtaining bounds for the number of
possible slopes which are independent of the particular manifold
we are studying. It turns out that only the genus of the boundary of
$M$ is relevant.

We define the genus of $\partial M$ be the sum of the
genus of the components of $\partial M$.

\proclaim{Theorem 2}
There is a function $n(g,g_\partial)$
such that there are at most $n(g,g_\partial)$ $\partial$-slopes
for essential surfaces of genus $g$ in a $\partial$-irreducible,
acylindrical, atoroi-
dal 3-manifold whose boundary has genus equal to
$g_\partial$.
\endproclaim 

We can also obtain bounds on the number of boundary slopes
in terms of hyperbolic geometry.

{\bf Definition.} Let $\Cal M (V)$ be the set of all hyperbolic
3-manifolds of totally geodesic boundary and with volume bounded
above by $V>0$. 

\proclaim{Theorem 3} There is a function $n_1(g, V)$ such that
there are at most $n_1(g,V)$
$\partial$-slopes for essential surfaces of genus $g$
in a 3-manifold $ M\in \Cal M(V)$. \endproclaim

Theorem 3 follows from either Theorem 1 or Theorem 2,
the fact that all maximum torus cusps have volumes $>C>0$ [Ad], 
and the following

\proclaim{Theorem 4}
There is an integer $g^*>0$ and a number $L >0$ such that if $M\in \Cal M(V)$,
then 

\noindent
(1) the genus of $\partial M$ is at most $g^*$. 

\noindent
(2) the length of any closed geodesic on $\partial M$ is at least $L$. \endproclaim

\noindent
{\bf Remark on Theorem 4.} Theorem 4 can be restated as follows:
For hyperbolic 3-manifolds with totally geodesic boundary and bounded volume,
the areas of their boundaries have an upper bound, and the lengths
of simple closed geodesics on their boundary have a lower bound.
Neither of those two assertions is true in dimension 2. Surfaces
of given area can have geodesic boundaries of any length. 

\demo{Proof of Theorem 1}
If $M$ has any 2-sphere boundary components,
we can fill them in with balls without changing the number of boundary slopes.
Since any essential surface can be homotoped off of a splitting 2-sphere,
we can without loss of generality assume that $M$ is irreducible.
The number of boundary slopes of essential surfaces lying on a torus
boundary component of $M$ is finite by [HRW], so we restrict attention
to surfaces with boundary on a higher genus component of $\partial M$.
By Thurston's Geometrization
Theorem for Haken manifolds, $M$ admits a complete
hyperbolic structure of finite volume with totally geodesic boundary [T].
We assume that $M$ is equipped with such a hyperbolic structure.
The totally geodesic boundary components consist of the non-torus
boundary components of $M$.
Since $\partial M$ have only finitely many components,
to prove Theorem 1, we need only to show that for each component
of $\partial M$
there are finitely many boundary slopes of proper
essential surfaces of genus at most $g$.

Suppose $F$ is an incompressible, boundary incompressible proper
immersion with $\partial F$ consisting of $n$ copies of a slope $l$.
Let $DM$ be the double of $M$ along its totally geodesic boundary components.
$DM$ is Haken and atoroidal, and admits a hyperbolic structure obtained
by doubling that of $M$. The double $DF$ of $F$ is incompressible,
and therefore a theorem of Schoen-Yau and Sacks-Uhlenbeck shows
that there is a least area representative of its homotopy class,
denoted by $DF^*$ [SY]. The intersection of $DF^*$ with the incompressible
least area (in fact totally geodesic) surface $\partial M$ consists of
curves essential on both $DF^*$ and $\partial M$.
Since $F$ is boundary incompressible in $M$,
the intersection $F^*$ of $DF^*$
with $M$ is a surface homotopic (rel boundary) in $M$ to $F$.
Since $DM$ admits an isometry which is a reflection about $\partial M$,
$DF^*$ is perpendicular to $\partial M$. If not,
we could reflect $DF^* \cap M$ and get a homotopic surface with lower area.
So $F^*$ is properly homotopic to $F$, $F^*$ is perpendicular to
$\partial M$ and $\partial F^*$ is a (possibly multiply covered) geodesic. 

Choosing geodesic orthogonal coordinates near the geodesic boundary
of the surface $F^*$, we have (line 7 of p.374, [BM])
$$ ds^2=du^2 + J^2(u,v)dv^2,\qquad (J(u,v)>0 \text{ and } J(0,v)=1).\tag 1 $$
where the $u$-curves (those where $v=$ constant) are geodesics
perpendicular to the boundary and the $v$-curves lie on the boundary when $u=0$.

The geodesic curvature in $F^*$
of a curve $t\mapsto (u(t), v(t))$ is given by Formula 10.4.7.1 of [BM], 

$$
\frac 1{\sqrt{E{u'}^2+G{v'}^2}}(\frac {d\phi}{dt}+ \frac 1{2\sqrt{EG}}
(\frac{\partial G}{\partial u}v'- \frac{\partial E}{\partial v}u')),\tag 2 $$
where $\phi$ is the angle between the curve and the $u$-curves and
the metric on $F^*$ is given by
$$E du^2 + G dv^2 .$$

When we consider the $v$-curves, we have $u'=0$, $v'=1$, $\phi=\pi/2$, $E=1$
and $G=J^2$. Substituting into (2), the geodesic curvature for a
$v$-curve $\{ u = c \}$ oriented as the boundary of $\{ 0 \le u \le c \}$
is given by:
$$k_g=\frac {1}{J} \frac {\partial J}{\partial u}. \tag 3$$
Orienting the curve as the boundary of $\{ u \ge c \}$ changes the
sign and gives
$$k_g=-\frac {1}{J} \frac {\partial J}{\partial u}. \tag 3'$$ 

The Gaussian curvature of the surface is Formula 10.5.3.3 of [BM] 

$$K=-\frac 1J \frac {\partial^2 J}{\partial u^2}.\tag 4$$ 

A direct computation shows that $k_g$ satisfies the following equation 

$$ \frac {\partial k_g}{\partial u}=K+k_g^2.\tag 5$$ 

Since $M$ is of constant curvature $-1$, we have 

$$K=k_1k_2-1$$
by Gauss's Formula (p.179 [Sp]), where $k_1$ and $k_2$ are the principle
curvatures. Since $F$ is a minimal surface, we have $k_1k_2\le 0$,
and hence $K\le -1$. Then by (4) it follows that
$$ \frac {\partial^2 J}{\partial u^2}\ge J.\tag 6$$ 

Fixing $v=v_0$, by (6) we have
$$\frac {\partial J}{\partial u}(u,v_0)
=\frac {\partial J}{\partial u}(u,v_0)
-\frac {\partial J}{\partial u}(0,v_0)$$
$$=\int^u_0\frac {\partial^2 J}{\partial u^2}(s,v_0)ds\ge \int^u_0 J(s,v_0)ds
\ge 0.\tag 7$$ 

(1), (3') and (7) imply that $k_g<0$, if $u>0$. 

Now consider the function
$$h(u)=\frac{e^{-u}-e^u}{e^{-u}+e^u}.\tag 8$$

which is the solution to the differential equation

$$\frac{dh}{du}=-1+h^2\tag 9$$

with the initial condition $h(0)=0$. Note that $h(u)<0$ when $u>0$ and
that the function $k_g-h$ satisfies the differential inequality 

$$\frac{d(k_g-h)}{du}=K+k_g^2+1-h^2$$
$$\le k_g^2-h^2=(k_g-h)(k_g+h)\tag 10$$
by (5) and (9). We want to show that $k_g-h\le 0$.

Suppose on the contrary that on an interval $[0, U]$,
$k_g-h$ is somewhere positive. Pick a $u_0\in [0,U]$,
such that $k_g-h$ takes its positive maximum at $u_0$.
We know $u_0\ne 0$ since $\partial F^*$ is a geodesic
and $k_g=h=0$ at $0$. Then
$$ \displaystyle \frac{d(k_g-h)}{du}$$
is zero if $u_0\in (0,u)$, and is $\ge 0$ if $u_0=U$.
Hence
$$ \displaystyle \frac{d(k_g-h)}{du} \ge 0$$ at $u_0$.
Since both $k_g$ and $h$ are negative at $u_0$,
we have $$(k_g-h)(k_g+h)< 0.$$ This contradicts (10), and so $k_g\le h$. 

For $t>0$, let $N_t(\partial M)$
be the subset of $M$ with distance $\le t$ from the boundary.
There is a $b>0$ such that when $t<b$ then $N_t(\partial M)$
is a collar of $\partial M$.

Choose $U<b$ in the above and let $N_U(\partial F^*)$ be the neighborhood of
$\partial F^*$ with $u$ coordinates at most $U$. Clearly
$N_U(\partial F^*)\subset N_b(\partial M)$. Since $N_b(\partial M)$
is a collar of $\partial M$ and the surface $ F^*$ is
$\partial$-incompressible, it follows that $N_U(\partial F^*)$
is a collar of $\partial F^*$. Letting
$$F_U=\overline{F^*-N_U(\partial F^*)},$$
each component of $\partial F_U$ is in the same homotopy class in $M$
as the slope $l$ and $$\# \partial F_U=\# \partial F^* =\# \partial F=n.$$ 

By Gauss-Bonnet, we have that

$$\int_{F_U} KdA +\int_{\partial F_U} k_gds =2\pi(\chi(F))=2\pi(2-2g-n).$$

Let $d$ be the length of the geodesic in the homotopy class of the slope $l$.
Then the length of each component of $\partial F_U$ is larger than $d$.
Since $K\le -1$ and $k_g\le h< 0$ at $U$, we have

$$nhd\ge 2\pi(2-2g-n).$$
Then we have

$$ d \le\frac{2\pi(2g+n-2)}{-hn}\le \frac{2\pi(2g+1)}{-h}.\tag 11$$

Since $g$ is given and $h = h(U) <0$, $d$ is bounded above.
There are only finitely many homotopy class of
essential closed curves in $\partial M$ containing elements of length less
than a given constant. Therefore for any fixed $g$
there are only finitely many
$\partial$-slopes for immersed incompressible,
boundary incompressible surfaces of genus $g$. \qed \enddemo 

\demo{Proof of Theorem 2}
In the proof of Theorem 1,
if we only consider 3-manifolds whose totally geodesic boundary has a collar
of width bounded below by $U_*$, so that $U > U_*$, then
$h=h(U)\le h(U_*)=-\tanh U_*\le 0$. Moreover if we consider only boundary
slopes of length at least $L>0$, then by (11) we have
$$L \le d \le \frac{2\pi(2g+1)}{\tanh U_*}.\tag 12$$

Let $A(R)$ be the area of $D(R)$, the hyperbolic disc with radius $R$.
Let $\Gamma_{L}$ be any lattice on the hyperbolic plane such that the
distance of any two vertices has distance at least $L$. Then the number
of vertices of $\Gamma_{L}$ in $\displaystyle D(\frac {2\pi(2g+1)}{\tan U_*})$
is at most
$$ \displaystyle n(g, U_*, L)=\frac{A(\frac{2\pi(2g+1)}{\tan U_*}+L)}{A(L)}.\tag 13$$
It follows that the number of boundary slopes for proper essential surfaces
of genus at most $g$ is bounded by $n(g, U_*, L)$.

To show the existence of the function $n(g,g_\partial)$
we need to establish in the proof of the previous theorem: 

\noindent
1. A lower bound $U_*$
to the width of a collar around $\partial M$ for any hyperbolic metric
on a manifold $M$ in which $\partial M$ is totally geodesic of genus
$\le g_\partial$. 

\noindent
2. Given $L >0$, an upper bound on the number of curves of length
$ \le L$ lieing in a collar of $\partial M$.
This bound should depend only on the genus of $\partial M$,
and not on its geometry. 

The existence of the first type of bound was established by Kojima and Miyamoto
[KM], and by Basmajian [Ba].
On the boundary of $M$, the second type of bound is a consequence of
the Margulis Lemma, or of its two dimensional version known as the
``collar lemma'' ([Bu] and also Theorem 2.18 of [Mu]).
We actually use a bound that holds in a collar neighborhood of
$\partial M$ in Theorem 1.
However the projection from a collar of the boundary of a
hyperbolic manifold with totally geodesic boundary to the boundary
is length decreasing,
so it suffices to consider curves lieing on the boundary. 

More precisely, let $S(x)=\sinh^{-1}(1/\sinh (x/2))$.
For a given simple closed geodesic $c$ with length $d_c$
on a hyperbolic surface, let $N(c)=\{x: d(x, c)\le S(d_c)\}$.
Then the collar lemma states that $N(c)$ is a collar.
Moreover if $c_1$ and $c_2$ are disjoint simple closed geodesics,
then $N(c_1)$ and $N(c_2)$ are disjoint.
There is an $L$ such that if $d\le L$, then $S(d) > d/2$
and $d> \sinh (d/2)$; for example, we can choose $L=1.75$;
then $S(d)> S(L)> 0.887> 0.85=L/2\ge d/2$ and
$1.76 \le L/\sinh (L/2)\le d/\sinh(d/2)$.
Then any two simple closed geodesics of length $\le L$ are disjoint.
Moreover the area of $N(c)$ is
$$
2d_c \sinh(S(d_c))= 2d _c/\sinh (d_c/2)\ge 2 d_c/d_c=2. $$
Hence the number of simple closed geodesics of length at most
$L$ is bounded above by 

$$2\pi(2g(F)-2)/2= 2\pi(g(F)-1).\tag 14$$ 

where $g(F)$ is the genus of $F$.

For simplicity, we first assume that $\partial M$ is connected.
By (13) and (14) we have

$$n(g, g_\partial)=\frac{A(\frac{2\pi(2g+1)}{\tanh U^*}+L)}{A(L)} + 2\pi(g_\partial-1).\tag 15$$

By Lemma 3.1 of [Ba], we have the lower bound 

$$U^*= \frac 14 log \frac{g_\partial+1}{g_\partial-1}.\tag 16$$
Moreover $\displaystyle A(R)=\frac {4\pi} {1- \tanh^2 R/2}$.
We can get an explicit value for
$n(g,g_\partial)$ by plugging in these functions,
though this does not appear to give sharp values. 

In general suppose $\partial M$ consists of $k$ torus components and $l$
components of genus $g_i>1$, $i=1,...,l$.
Then $g_\partial=\sum_{i=1}^l g_i+ k$ and
there are at most $\sum_{i=1}^l n(g, g_i)+k N(g)$
boundary slopes for proper essential surfaces of genus $g$,
where $N(g)$ is the uniform bound for the number of
boundary slopes of proper essential surface of genus $g$
on a torus boundary component given in [HRW], and $n(g, g_i)$ is given by (15).
One can verify that $\sum_{i=1}^l n(g,g_i)+N(g)\le n(g, g_\partial)$.
\qed
\enddemo

\demo{Proof of Theorem 4}
Pick any infinite sequence of totally geodesic hyperbolic
3-manifolds $\{M_n\}$ in $\Cal M(V)$. Consider the sequence
$\{D(M_n)\}$, where $D(M_n)$ is the double of each $M_n$.
Then the volume of the closed hyperbolic 3-manifold $D(M_n)$
is bounded by $2V$. By passing to a subsequence,
we can assume that $D(M_n)$ has a Gromov limit $M^*$. It is known that 

\noindent
(a) $M^*$ is a complete hyperbolic 3-manifold of finite volume,
which can be viewed as the complement of a hyperbolic link $L$ in a
closed 3-manifold, and each $D(M_n)$ is obtained by a Dehn surgery on $M^*$. 

\noindent
(b) Since each $D(M_n)$ admits a reflection $r_n$ (isometry) about
its geodesic boundary, so does $M^*$. Hence $M^*=D(M_\infty)$,
where $M_\infty$ is a hyperbolic 3-manifold with totally geodesic boundary.
Let $r_\infty$ be the reflection of $D(M_\infty)$ about $\partial M_\infty$.
We have not claimed as yet that there is no cusp at $\partial M_\infty$. 

\noindent
(c) Let $TH_\epsilon(P)$ be the $\epsilon$ thick part of $P$
for any hyperbolic 3-manifold $P$. Then for any $\epsilon>0$ and
$1-\epsilon \le k\le 1$, there is an integer $N$ such that for
$n>N$ there is a homeomorphism
$h_n: TH_\epsilon (D(M))\to TH_\epsilon(D(M_n))$ which is a $k$-quasi-isometry. Moreover $h_n$ can be chosen to commute with the reflections. 

For the result on the Gromov limit of closed hyperbolic
3-manifolds of bounded volume, see Chapter 6 of [T1], or Chapter E of [BP].
For the fact about reflections, one can argue as follow:
As in the case of closed hyperbolic 3-manifolds, any sequence of hyperbolic
3-manifolds with totally geodesic boundary and bounded volume $V$ has
a subsequence with Gromov limit $M_\infty$, which is a complete
hyperbolic 3-manifold with totally geodesic boundary.
Then the double $D(M_\infty)$ will be the limit of the doubles.

Suppose there is a torus in $\partial M$, which must result
in cusps in $\partial M_\infty$. Let the torus $T$ be the boundary
component of $TH_\epsilon (M_\infty)$ corresponding to the cusp $C$,
and let $c$ be a component of $T\cap \partial M$. Then $h_n(T)\subset D(M_n)$
is invariant under the reflection $r_n$ about $\partial M_n$,
and it follows that $h_n(c)$ is a meridian of the Dehn filling solid torus
on $h_n(T)$, and therefore $h_n(c) \subset \partial M_n$ is a trivial loop.
However each cusp in $\partial M_\infty$ can only be a limit of
essential loops, and this is a contradiction. Hence $\partial M_\infty$
contains no cusps.

Since $\partial M_\infty$ contains no cusps,
for small $\epsilon$, $\partial M_\infty$ is contained in the
interior of the compact manifold $TH_\epsilon D(M_\infty)$.
Moreover as the fixed point set of the reflection
$r_\infty | TH_\epsilon (D(M_\infty))$, $\partial M_\infty$ is compact,
therefore it is closed. Since $\partial M_n$ converges to $\partial M_\infty$
in the limit, it follows that 

\noindent
(1') the genus of $\partial M_n$ is stable when $n$ is large enough, 

\noindent
(2') the length of the shortest simple closed geodesic
on $\partial M_n$ cannot converge to zero
(otherwise there will be a cusp in $M_\infty$). 

Now suppose (1) of Theorem 4 is not true. Then we can find a
sequence $\{M_n\}$ in $\Cal M(V)$ such that the genus of
$\partial M_n$ is $>n$. The genus of any subsequence must also
tend to infinity, which contradicts (1'); hence (1) of Theorem 4 is true.
Similarly, if (2) of Theorem 4 is not true,
then we can find a sequence $\{M_n\}$ in $\Cal M(V)$ such
that the length of the shortest geodesic of $\partial M_n$ is $<1/n$,
which contradicts (2'). This finishes the proof of Theorem 4. \qed \enddemo

\vskip 0.5 true cm

{\bf References.}

\vskip 0.5 true cm
[Ad] C. Adams, Volumes of n-cusped hyperbolic 3-manifolds,
J. London Math. Soc. 1988, 38, 2, 555-565. 

[Agol] I. Agol, Topology of Hyperbolic 3-manifolds, Ph.D. thesis, UCSD, 1998. 

[Ba] A. Basmajian, Tubular neighborhoods of totally geodesic
hypersurfaces in hyperbolic manifolds, Invent. Math. 117 (1994), 207-225. 

[BP] R. Benedetti and C. Petronio, Lectures on Hyperbolic Geometry,
Universitext, Springer-Verlag, (1991). 

[Bu] P. Buser,
The collar theorem and examples. Manuscripta Math. 25 (1978), 349--357.

[BG] M. Berger and B. Gostiaux, Differential Geometry:
Manifolds, curves and surfaces, GTM 115, Springer Verlag, Berlin, New York. 

[Go] C. Gordon, Dehn surgery on knots,
Proceedings of the International Congress of Mathematicians,
Vol. I, II (Kyoto, 1990), 631--642, Math. Soc. Japan, Tokyo, 1991. 

[HRW] J. Hass, H. Rubinstein and S.C.Wang,
Immersed surfaces in 3-manifo-
lds, preprint (1998).

[Ha] A. Hatcher, On the boundary curves of incompressible surfaces,
Pacific J. Math. 99 (1982), 373-377. 

[KM] S. Kojima and Y. Miyamoto,
The smallest hyperbolic $3$-manifolds with totally geodesic boundary.
J. Differential Geom. 34 (1991), 175--192. 

[Lu] J. Luecke, Dehn surgery on knots in $S^3$,
Proc. ICM Vol 2 (Zurich, 1994), 585-594. 

[Mu] C. McMullen, Complex dynamics and renormalization,
Ann. Math. Study, No, 135, PUP, 1994. 

[Oe] U. Oertel, Boundaries of $\pi_1$-injective surfaces,
Topology Appl. 78, (1997), 215-234. 

[SW] M. Scharlemann and Y. Wu,
Hyperbolic manifolds and degenerating handle additions,
J. Aust. Math. Soc. (Series A) 55 (1993), 72-89. 

[Sh] P. Shalen, Representations of 3-manifold groups and
its application to topology, Proc. ICM Berkeley, (1986), 607-614. 

[Sp] M. Spivak, A comprehensive Introduction to Differential
Geometry, Vol. 4, Publish or Perish, Inc. Berkeley 1979. 

[SY] R. Schoen and S.T. Yau,
Existence of incompressible minimal surfaces and the
topology of three-dimensional manifolds with nonnegative
scalar curvature, Ann. of Math. (2) 110 (1979), 127-142. 

[T] W. Thurston, Three dimensional manifolds, Kleinian
groups and hyperbolic geometry, Bull. AMS, Vol. 6, (1982) 357-388. 

[T1] W. Thurston, Geometry and Topology of 3-manifolds,
Princeton University Lecture Notes, 1978.

\bye